\documentclass[reqno,final]{amsart}
\usepackage{natbib}  
\usepackage{fancyhdr} 
\usepackage{color} 
\usepackage{hyperref} 
\usepackage{graphicx} 
\usepackage{enumerate}

\definecolor{aleacolor}{rgb}{0.16,0.59,0.78}

\hypersetup{
breaklinks,
colorlinks=true,
linkcolor=aleacolor,
urlcolor=aleacolor,
citecolor=aleacolor}

\renewcommand{\headheight}{11pt}
\renewcommand{\cite}{\citet}

\theoremstyle{plain}
\newtheorem{theorem}{Theorem}[section]                                          
\newtheorem{prop}[theorem]{Proposition}                          
\newtheorem{lemma}[theorem]{Lemma}

\theoremstyle{definition}

\theoremstyle{remark}

\makeatletter \@addtoreset{equation}{section} \makeatother
\renewcommand\theequation{\thesection.\arabic{equation}}
\renewcommand\thefigure{\thesection.\arabic{figure}}
\renewcommand\thetable{\thesection.\arabic{table}}

\newcommand{\aleaIndex}[1]{\href{http://alea.impa.br/english/index_v#1.htm}{\bf #1}}
\eheader{Alea}{\aleaIndex{7}}{2010}{19}{40}

\elogo{\parbox[c]{3cm}{\includegraphics[width=3cm]{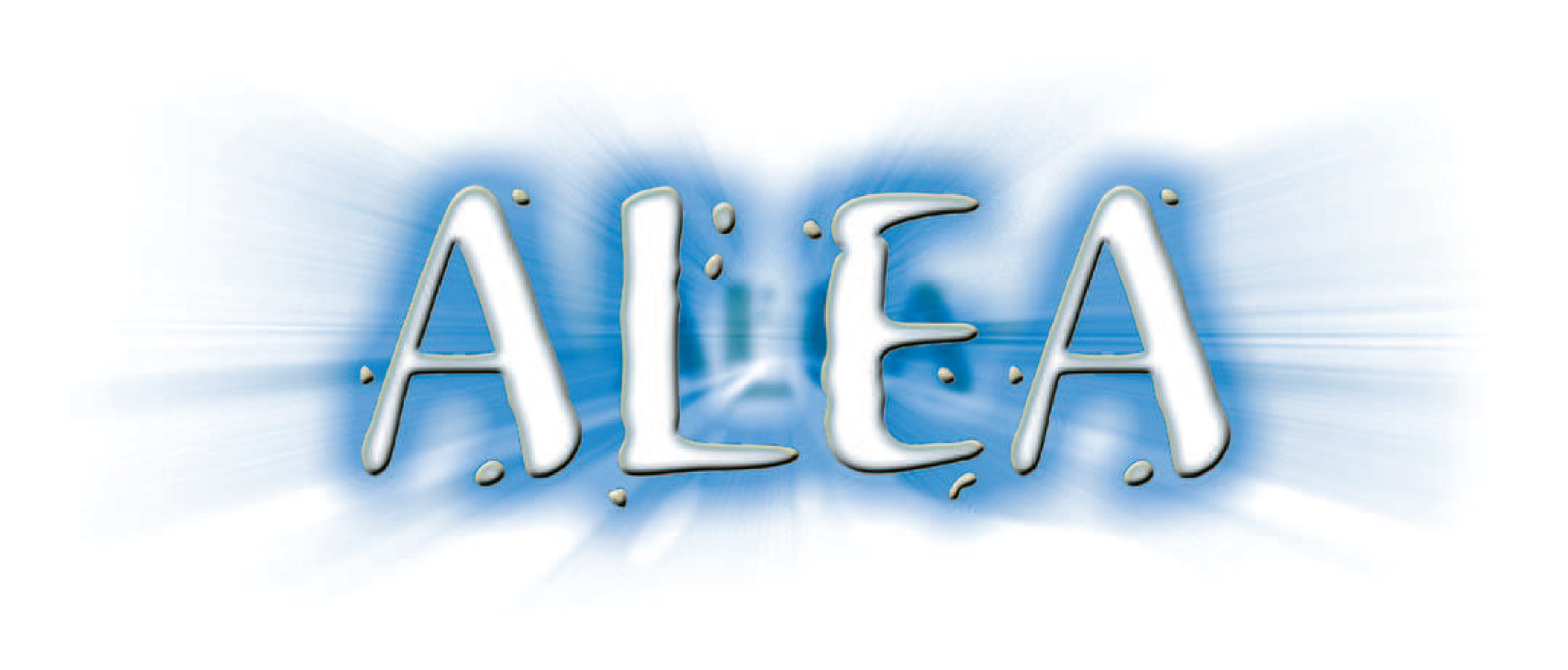}}}

\usepackage{epsfig,amsmath,amsfonts,latexsym,subfigure}

\newcommand{\E}{{\mathbf{E}}}

\newcommand{\Z}{{\mathbb Z}}

\newcommand{\Pb}{{\mathbb P}}
\newcommand{\Pf}{{\mathbf{P}}}
\def\Cox{\hfill \Box}
\def\eps{\varepsilon}

\begin{document}

\date{June 5, 2009; accepted March 17, 2010}

\author{Maria Deijfen}
\address{Stockholm University}
\urladdr{www.math.su.se/$\sim$mia}
\email{mia@math.su.se}

\author{Olle H\"{a}ggstr\"{o}m}
\address{Chalmers University of Technology}
\urladdr{www.cs.chalmers.se/$\sim$olleh}
\email{olleh@chalmers.se}

\title[Biased random walk in translation invariant percolation]{On the speed of biased random walk in translation invariant percolation}

\begin{abstract}

\noindent For biased random walk on the infinite cluster in
supercritical i.i.d.\ percolation on $\Z^2$, where the bias of the
walk is quantified by a parameter $\beta>1$, it has been conjectured
(and partly proved) that there exists a critical value $\beta_c>1$
such that the walk has positive speed when $\beta<\beta_c$ and speed
zero when $\beta>\beta_c$. In this paper, biased random walk on the
infinite cluster of a certain translation invariant
percolation process on $\Z^2$
is considered. The example is shown to exhibit the opposite behavior
to what is expected for i.i.d.\ percolation, in the sense that it has a
critical value $\beta_c$
such that, for $\beta<\beta_c$, the random walk has speed zero, while,
for $\beta>\beta_c$, the speed is positive. Hence the monotonicity
in $\beta$ that is part of the
conjecture for i.i.d.\ percolation cannot be extended to general translation
invariant percolation processes.
\end{abstract}

\maketitle

\section{Introduction}

This paper is concerned with biased random walk on infinite percolation
clusters on the square lattice, whose vertex set is $\Z^2$ and whose
edge set consists of pairs of vertices at Euclidean distance $1$ from
each other; with a slight abuse of notation we write $\Z^2$ for this
lattice. Let there be two possible states for each edge
in $e \in E$: open or closed. In general, a percolation model is a way
of deciding which edges are to be open. In standard i.i.d.\ bond
percolation with parameter $p\in[0,1]$, each edge is independently
open with probability $p$. The resulting configuration
will almost surely contain an infinite open cluster if and
only if $p > 1/2$; see 
\cite{MR1707339} for this and other basics on percolation theory.
When the origin belongs to the infinite cluster, we can define a random walk, starting at the origin, as follows: Let
$Z_t=(X_t,Y_t)$ denote the position of the random walk at time
$t\in\mathbb{N}$ (we apologize to sensitive readers for using the
letter $t$ for a discrete time parameter; however, the integer indices
$i,j,k,l,m,n,\ldots$ will be needed for other purposes later on).
Write $\Lambda_t$ for the set of neighbors of
$Z_t$ in the infinite cluster and define $l_t=|\Lambda_t|$. Also,
fix $\beta>1$. If $(X_t+1,Y_t)\in\Lambda_t$, then $Z_{t+1}=(X_t+1,Y_t)$
with probability $\beta(\beta+l_t-1)^{-1}$ and $Z_{t+1}$ equals any
other given vertex in $\Lambda_t$ with probability $(\beta+l_t-1)^{-1}$.
If $(X_t+1,Y_t)\not\in \Lambda_t$, then $Z_{t+1}$ is chosen uniformly from
$\Lambda_t$, that is, $Z_{t+1}$ equals any given vertex in $\Lambda_t$
with probability $l_t^{-1}$. In case $\Lambda_t= \emptyset$, the walk stays
put, i.e., $Z_{t+1}=Z_t$.

This model was introduced by \cite{0022-3719-16-8-014} and describes a
random walk with drift towards the right, the strength of the drift
being quantified by the parameter $\beta$. (Note that zero drift
corresponds to $\beta=1$.) The asymptotic speed, or simply the speed, is
defined as $\lim_{t\to\infty}X_t/t$ (provided the limit exists).
In \cite{0022-3719-16-8-014}, it
is conjectured that there is a critical drift $\beta_c=\beta_c(p)>1$
such that the walk has positive speed for $\beta<\beta_c$ and speed
zero for $\beta>\beta_c$. Intuitively, if the drift is large, the
walk will tend to get stuck in ``dead ends'' of the percolation cluster,
while, if the drift is weaker, it will be able to quickly backtrack and
get out of the dead ends. The conjecture from \cite{0022-3719-16-8-014}
was partly confirmed in two simultaneous and independent papers by
\cite{MR2004982} and \cite{MR1990055}, respectively, where it is proved that
there are $\beta_l$ and $\beta_u$, with $1<\beta_l\leq \beta_u$,
such that the walk has positive speed for $\beta<\beta_l$ and speed
zero for $\beta>\beta_u$. (Sznitman in fact obtained the same result in
arbitrary dimension $d \geq 2$.)
What remains here is to show that one can take
$\beta_l = \beta_c$. \cite{MR2568279} demonstrated
the same critical phenomenon with $\beta_l = \beta_c$ for a certain dependent
percolation model on the lattice sometimes known as the infinite ladder.
One might ask whether the monotonicity property suggested by the
Barma--Dhar conjecture (namely that zero speed at a given $\beta>1$ implies
the same thing at all larger values of $\beta$) should be extended to a
wider class of percolation processes, such as those that are translation
invariant. Our main result, Theorem \ref{th:main} below, shows that the
answer is no.

More precisely, what we do in the present paper is as follows. We will
construct a translation invariant percolation
process on $\Z^2$ for which the above random walk dynamics give rise to a
process which has speed zero when $\beta$ is small and positive
speed when $\beta$ is large. For a translation invariant
probability measure $\Psi$
on $\{{\rm open,closed}\}^E$ defining the percolation process, with the
property that the existence of an infinite cluster has probability $1$,
write $\Pb_{\Psi,\beta}$ for the joint law of the percolation configuration
and the random walk $\{Z_t\}_{t\geq 0}$. Furthermore, write $\{0 \leftrightarrow \infty\}$
for the event that the origin belongs to an infinite
open cluster of the percolation configuration.
We will prove the following:

\begin{theorem}\label{th:main}
For each $\gamma>1$ there exists a translation invariant probability
measure $\Psi=\Psi(\gamma)$ on $\{{\rm open,closed}\}^E$ such that,
for any $\beta >1$,
$$
\lim_{t\to\infty}\frac{X_t}{t}= \left\{
\begin{array}{lll}
0 & \Pb_{\Psi,\beta}\textrm{-a.s. on the event }
\{0 \leftrightarrow \infty\}& \mbox{if }\beta<\beta_c \\
\frac{\beta-1}{\beta+1} & \Pb_{\Psi,\beta}\textrm{-a.s. on the event }
\{0 \leftrightarrow \infty\} &
\mbox{if }\beta>\beta_c
\end{array} \right.
$$
with $\beta_c= \gamma$.
\end{theorem}

\noindent To work out the speed at criticality is probably doable with some more work, but might not be so important, since we suspect that you can get either answer (zero speed or full speed) by further fine-tuning of our model.

The rest of the paper is organized as follows.
A percolation process with the property described in
Theorem \ref{th:main} is constructed in detail in Section
\ref{sect:main_construction}. First, however
we illustrate one of the main ideas by describing
in Section \ref{sect:warm-up} a simpler
(but not translation invariant) percolation process on which biased
random walk behaves as in the theorem; the key concept here is the ``trap''
structure in Figure 1 below. These traps appear also in the main construction
in Section \ref{sect:main_construction}. Section \ref{sect:without_traps}
concerns the main construction minus the traps, where
the asymptotic speed is shown to equal $\frac{\beta-1}{\beta+1}$ for any
$\beta>1$. In Section \ref{sect:positive_speed}
we show that including the traps as in the main
construction makes no difference to the asymptotic speed as long as
$\beta>\beta_c$, thus establishing
the second half of Theorem \ref{th:main}. Finally, in Section
\ref{sect:zero_speed}, we consider
the case $\beta<\beta_c$, and show that the traps slow down the speed to zero,
thereby proving the first half of the theorem.

\section{A warm-up construction} \label{sect:warm-up}

Consider first a configuration of open edges outlined in Figure 1(a): an
infinite open path starting at the origin and going off straight along
the positive $x$-axis, with so called
\emph{traps} attached to it.
Each trap consists of an \emph{entrance} and a \emph{core}, the core
being located one floor above the entrance; see Figure 1(b). The
length (meaning the number of edges) of the entrance and the core of the
$n$'th trap are denoted by $e_n$ and $c_n$ respectively. We furthermore write
$d_n$ for the $x$-coordinate where the entrance of the
$n$'th trap begins, and define $\Delta d_n= d_n-d_{n-1}$.
For each $n$, we need to have
\begin{equation} \label{eq:avoid_congestion_1}
e_n<\Delta d_n
\end{equation}
and
\begin{equation}  \label{eq:avoid_congestion_2}
e_n + c_{n-1}-e_{n-1} < \Delta d_n
\end{equation}
in order for the traps not to overlap. The vertex $(d_n, 0)$ is called
the \emph{anchor} of the trap.

\begin{figure}
\centering\mbox{\subfigure[Line with traps attached to it.]{\epsfig{file=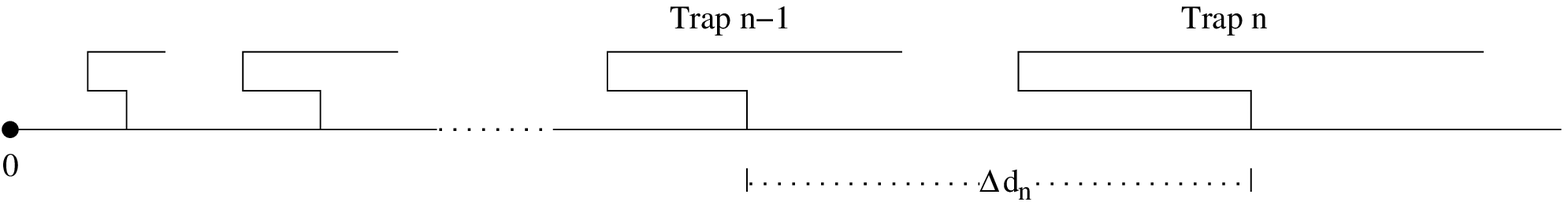,height=1.60cm}}}\par\par
\mbox{\subfigure[The $n$th trap = trap of size $n$.]{\epsfig{file=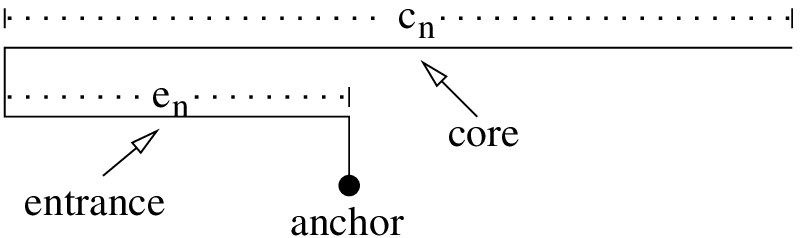,height=1.8cm}}}\caption{Schematic picture of configuration with line and traps.}
\end{figure}

It is readily checked that a random walk with drift $\beta>1$ along
the infinite line without traps has positive speed equal to
$\frac{\beta-1}{\beta+1}$. In particular, the walk is transient, and
adding traps cannot change this, due to Rayleigh's monotonicity principle;
see, e.g.,\cite{MR920811}.

To obtain a configuration on which the walk
has speed zero for small $\beta$ and positive speed for large $\beta$
we will choose the sequences $\{e_n\}$, $\{c_n\}$ and $\{d_n\}$
so that, if $\beta$ is small, the walk will enter the core of infinitely
many traps, which will cause a delay severe enough to bring the speed
down to zero, while, if $\beta$ is large, the walk will enter the core
of only finitely many traps, causing a delay that is negligible in the limit.
The details are as follows.

Fix $\alpha>0$. It will turn out that choosing
\begin{equation}  \label{eq:parameters_in_warmup_example}
d_n=n^3, \, \, \, e_n= \lceil \alpha \log n \rceil \,\,\, \mbox{and} \,\,\,
c_n=n
\end{equation}
for all $n$ sufficiently large (where $\lceil \cdot \rceil$ denotes
rounding up to the nearest integer) yields an asymptotic speed which is zero
for $\beta<\beta_c$ and strictly
positive for $\beta>\beta_c$, with $\beta_c=e^{1/\alpha}$.
The reason we say ``for all $n$ sufficiently large'' is that the
values of $e_n$ and $c_n$ may need to be lowered compared
to (\ref{eq:parameters_in_warmup_example}) for small $n$ in order to
satisfy (\ref{eq:avoid_congestion_1}) and (\ref{eq:avoid_congestion_2});
by transience of the random walk, modifying a finite number of traps
cannot change the asymptotic speed.
\begin{prop}  \label{prop:warm-up_construction}
The infinite path on the positive $x$-axis, decorated with traps
with parameters satisfying (\ref{eq:parameters_in_warmup_example}) for large
$n$, yields the following almost sure behavior for random walk:
\[
\lim_{t\to\infty}\frac{X_t}{t}= \left\{
\begin{array}{ll}
0 & \mbox{if }\beta<\beta_c \\
\frac{\beta-1}{\beta+1} &
\mbox{if }\beta>\beta_c
\end{array} \right.
\]
with $\beta_c= e^{1/\alpha}$.
\end{prop}
Before proving
this result, let us explain why it does not immediately imply our
desired Theorem \ref{th:main}. The reason, of course, is that the construction
in the present section is not translation invariant. Now, there is a standard
way of turning a non-translation invariant construction into a translation
invariant one, namely via random translation of the original construction.
In this case, the natural thing to do would be the following.
\begin{enumerate}[(a)]
\item
Make copies of the original configuration shifted $k$ steps vertically, for
$k=3,6,9,\ldots$, and for $k=-3, -6, -9,\ldots$.
\item
Shift the configuration resulting from (a) $K$ steps vertically,
where $K$ is chosen according to uniform distribution on
$\{0,1,2\}$, thus making the model
invariant under vertical translation.
\item
Shift the configuration resulting from (b) $L$ steps horizontally to
the left, where $L$ is chosen according to uniform distribution on
$\{1, \ldots, l\}$, and then
consider the limit as $l \rightarrow \infty$, thus making the model
invariant also under horizontal translation.
\end{enumerate}
The problem with this approach is that since
$\lim_{n \rightarrow \infty}\frac{d_n}{n}=\infty$
in (\ref{eq:parameters_in_warmup_example}), the density of trap
entrances goes to zero, and they will disappear on us in step (c). This
disappearance can be avoided by a more elaborate fractal-like construction
of the percolation process described in Section \ref{sect:main_construction}.

For the proof of Proposition \ref{prop:warm-up_construction}, we will
need three lemmas concerning the behavior of random walk on traps; these
lemmas will become useful also later on when we analyse random walk on our
main construction, in Sections \ref{sect:positive_speed} and \ref{sect:zero_speed}.

For $n, i \geq 1$, let $T_{n,i}$ denote the time spent by the random walk
in trap number $n$ during its $i$'th visit to the trap, and let
\[
T^*_{n,i} \, = \left\{
\begin{array}{ll}
0 & \mbox{if the walk hits the trap's core during this visit} \\
T_{n,i} & \mbox{otherwise,}
\end{array} \right.
\]
with the convention that if the walk enters the trap exactly
$k$ times, then $T_{n,i}=T^*_{n,i}=0$ for all $i>k$.
The first lemma gives the probability, once a trap has been entered, of
reaching its core.
\begin{lemma}  \label{lem:hit_the_core}
Each time the random walk enters the $n$'th trap, it
has probability \\
$\frac{\beta-1}{\beta^{e_n+1}+\beta-2}$ of reaching
the core before exiting the trap, so that
\[
\Pf\left[ T_{n,i} > T^*_{n,i} \, | \, T_{n,i}>0 \right] \, =
\frac{\beta-1}{\beta^{e_n+1}+\beta-2}
\]
for any $n$ and $i$.
\end{lemma}
Once the walk enters the core, it has a fair chance of spending a very
long time (exponential in $c_n$) there.
The second lemma quantifies this.
\begin{lemma}   \label{lem:stay_in_core}
For each $n$ and $i$, we have
\begin{equation}  \label{eq:stay_in_core}
\Pf\left[ T_{n,i} \geq \beta^{c_n}  \, | \, T_{n,i}>0 \right] \, \geq
\frac{(\beta-1)^2}{2\beta(\beta^{e_n+1}+\beta-2)} \, .
\end{equation}
\end{lemma}
On the other hand,
provided the walk does {\em not} hit the core of the trap, its expected time
spent in the trap can be bounded uniformly in $n$. The third lemma makes
this precise. Let $\preceq$ denote stochastic
domination between random variables, i.e., $X \preceq X'$ means that
$\E[f(X)] \leq \E[f(X')]$ for any bounded and increasing $f$.
\begin{lemma}  \label{lem:exit_fast_if_core_not_hit}
We may define a positive random variable $T^*$ such that
\begin{description}
\item{\rm (a) } $T^*_{n,i} \preceq T^*$ for any $n$ and $i$, and
\item{\rm (b) } $\E[T^*] = \frac{2\beta-1}{\beta-1}$.
\end{description}
\end{lemma}
For the proofs of the lemmas (and also later on) it will be
useful to consider an electrical analysis of the random walk \`a la
\cite[Chapter 3]{MR920811}. Each edge $e$ of the network is assigned a
resistance $R(e)= \beta^{-x(e)}$ where $x(e)$ is the largest $x$-coordinate
amongst the two vertices incident to $e$. The rules of the random walk may
then be reformulated as saying that a random walker standing at a vertex
chooses among the incident edges with probabilities inversely proportional
to their resistances. For two vertices $v_1$ and $v_2$, let $R_{eff}(v_1, v_2)$
denote the effective resistance between $v_1$ and $v_2$ in the electrical representation, see \cite[ Section 3.4]{MR920811}. Since the percolation network in this case is a tree, there is always a unique self-avoiding path between $v_1$ and $v_2$, and $R_{eff}(v_1, v_2)$ is simply the sum
of the edge resistances along the path.

\medskip\noindent
{\bf Proof of Lemma \ref{lem:hit_the_core}.}
When the walk enters the trap, it will find itself at the vertex $(d_n,1)$.
From there, it will eventually reach either $(d_n,0)$ or $(d_n-e_n, 2)$. The
probability that it hits the latter before the former equals
\begin{eqnarray}
\nonumber
\frac{R_{eff}((d_n,1), (d_n,0))}{R_{eff}((d_n,1), (d_n,0))
+ R_{eff} ((d_n,1), (d_n-e_n, 2))} & = &
\frac{\beta^{-d_n}}{\beta^{-d_n}+ \sum_{i=d_n-e_n}^{d_n} \beta^{-i}} \\
& = &
\nonumber
\frac{1}{1+ \sum_{i=0}^{e_n} \beta^{i}} \\
& = & \frac{\beta-1}{\beta^{e_n+1}+\beta-2} \, .
\label{eq:first_resistance_comparison}
\end{eqnarray}
$\Cox$

\medskip\noindent
{\bf Proof of Lemma \ref{lem:stay_in_core}.}
Once the random walk hits the core, i.e., once it reaches the
vertex $(d_n-e_n, 2)$, its conditional
probability of hitting the second-to-last
vertex $(d_n-e_n+c_n-1,2)$ of the core before going back to
$(d_n-e_n, 1)$ is
\begin{eqnarray}
\nonumber
& & \mbox{ } \hspace{-12mm}
\frac{R_{eff}((d_n-e_n, 2), (d_n-e_n, 1))}{R_{eff}((d_n-e_n, 2),
(d_n-e_n, 1)) + R_{eff}((d_n-e_n, 2), (d_n-e_n+c_n-1,2))} \\
\nonumber
& = & \frac{1}{1+ \sum_{i=1}^{c_n-1} \beta^{-i}} \\
& \geq & \frac{1}{1+ \sum_{i=1}^{\infty} \beta^-i}
\, = \, \frac{\beta-1}{\beta} \, .
\label{eq:second_resistance_comparison}
\end{eqnarray}
Combining this with Lemma \ref{lem:hit_the_core}, we thus have that
once the random walk enters the trap, it has probability at
least
\begin{equation}  \label{eq:2nd_cond_prob}
\frac{(\beta-1)^2}{\beta(\beta^{e_n+1}+\beta-2)}
\end{equation}
of reaching $(d_n-e_n+c_n-1,2)$ before exiting. A similar calculation
as in (\ref{eq:first_resistance_comparison}) and
(\ref{eq:second_resistance_comparison}) shows that once the walk has reached
$(d_n-e_n+c_n-1,2)$, it has probability
$\frac{\beta-1}{\beta^{c_n+1}-1}<\beta^{-c_n}$ of hitting $(d_n-e_n,1)$ before
$(d_n-e_n+c_n, 2)$. Hence, upon reaching $(d_n-e_n+c_n-1,2)$, the number
of visits to $(d_n-e_n+c_n, 2)$ before
reaching $(d_n-e_n,1)$ is geometric with
mean at least $\beta^{c_n}$. Thus, the number of such visits exceeds
$\frac{\beta^{c_n}}{2}$ with conditional probability at least $\frac{1}{2}$,
and multiplying by (\ref{eq:2nd_cond_prob}) yields the corresponding
unconditional probability (which is the desired right-hand side of
(\ref{eq:stay_in_core})).
Every visit to $(d_n-e_n+c_n, 2)$ is immediately followed by one
to $(d_n-e_n+c_n-1,2)$, so by counting also the latter we can replace the count
$\frac{\beta^{c_n}}{2}$ by simply $\beta^{c_n}$, and (\ref{eq:stay_in_core})
follows. $\Cox$

\medskip\noindent
{\bf Proof of Lemma \ref{lem:exit_fast_if_core_not_hit}.}
Write $(Z'_1, Z'_2, \ldots, Z'_{T_{n,i}})$ for the sequence of vertices visited
during the $i$'th visit to trap $n$. Also, write
$(Z''_1, Z''_2, \ldots, Z''_{T^*_{n,i}})$ for the thinned sequence
obtained by deleting all visits to the core of the trap, and
note that $(Z''_1, Z''_2, \ldots, Z''_{T^*_{n,i}})$ has the same distribution
that the original sequence would have had if the trap had had no core.

Imagine now a trap whose entrance is infinite, i.e., consists of an
infinite straight path going off to the left from the anchor and generate
a sequence $(W'_{\!1}, W'_{\!2}, \ldots, W'_{\! T^*})$ describing the positions
of the random walk during a single visit to this trap. Since the
walk has a drift to the right ($\beta>1$), we get that $T^*< \infty$ with
$\E[T^*]< \infty$; a simple calculation shows that
$\E[T^*]= \frac{2 \beta-1}{\beta-1}$. Let $(W''_{1}, W''_{2}, \ldots, W''_{T^*_n})$
denote the thinned sequence obtained by deleting from
$(W'_1, W'_2, \ldots, W'_{T^*})$ all visits to vertices more than $e_n$ steps
to the left of the anchor, and note, crucially, that
$(W''_1, W''_2, \ldots, W''_{T^*_n})$ has the same distribution as
$(Z''_1, Z''_2, \ldots, Z''_{T^*_{n,i}})$. Hence $T^*_n$ and $T^*_{n,i}$
are identically distributed, and since $T^*_n \leq T^*$ the proof is complete.
$\Cox$

\medskip\noindent
Equipped with Lemmas \ref{lem:hit_the_core}, \ref{lem:stay_in_core} and
\ref{lem:exit_fast_if_core_not_hit}, we are now in a position to prove
Proposition \ref{prop:warm-up_construction}.

\medskip\noindent
{\bf Proof of Proposition \ref{prop:warm-up_construction}.}
For any vertex $(k,0)$ on the positive $x$-axis, the resistance
to infinity is given by
\[
R_{eff}((k,0),\infty) \, = \, \sum_{j=k+1}^\infty \beta^{-j}
\, = \, \frac{\beta^{-k}}{\beta-1} \, .
\]
It is a standard fact from \cite{MR920811} that the escape
probability $p_{esc}((k,0))$ of the random walk from $(k,0)$ -- that is, the probability that the walk leaves $(k,0)$ and never returns -- is given by
\begin{eqnarray}  \label{eq:escape_probability}
p_{esc}((k,0)) & = & \left( R_{eff}((k,0), \infty)
\sum_{\{e: \, e \sim (k,0)\}} R(e)^{-1}   \right)^{-1}
\end{eqnarray}
where $e \sim (k,0)$ means that the edge $e$ is incident
to the vertex $(k,0)$. If $k= d_n$ for some $n$ (i.e., there is some trap
connecting to the $x$-axis at $(k,0)$), we get
\begin{eqnarray*}
p_{esc}((k,0)) & = &
\left( \frac{\beta^{-k}}{\beta-1}
(\beta^{k+1}+ \beta^k + \beta^k)  \right)^{-1} \\
& = & \frac{\beta-1}{\beta+2} \, .
\end{eqnarray*}
For such $k$, the probability that the walk immediately takes a
step into the trap is $\frac{1}{\beta+2}$. Hence, the probability that
it ever takes a step into the trap before escaping to $\infty$ is
\[
\frac{\frac{1}{\beta+2}}{\frac{1}{\beta+2}+\frac{\beta-1}{\beta+2}}
\, = \, \beta^{-1} \, .
\]
It follows that
\begin{eqnarray}
\nonumber
& & \mbox{the number of visits to the trap is geometrically} \\
\label{eq:geometric_distr}
& & \mbox{distributed with mean $(\beta-1)^{-1}$.}
\end{eqnarray}
Define $T_{n, tot}$ as the total time $\sum_{i=1}^\infty T_{n,i}$ spent
in the $n$'th trap, and analogously
$T^*_{n, tot} =\sum_{i=1}^\infty T^*_{n,i}$.
Combining (\ref{eq:geometric_distr}) with
Lemma \ref{lem:hit_the_core} yields
\[
\Pf[T_{n, i} > T^*_{n, i}] \, = \,
\beta^{-i}\frac{\beta-1}{\beta^{e_n+1}+\beta -2}  \, .
\]
Summing over $i$ gives
that the probability $\Pf[T_{n, tot} > T^*_{n, tot}]$ of ever
hitting the core of the $n$'th trap satisfies
\begin{eqnarray}
\nonumber
\Pf [T_{n, tot} > T^*_{n, tot}] & \leq & \sum_{i=1}^\infty
\Pf [T_{n, i} > T^*_{n, i}] \\
\nonumber
& \leq & \frac{\beta-1}{\beta^{e_n+1}+\beta -2} \sum_{i=1}^\infty \beta^{-i} \\
& = & \frac{1}{\beta^{e_n+1}+\beta -2} \, .
\label{eq:total_time_in_core}
\end{eqnarray}

Consider first the case $\beta>\beta_c$, where we wish to show that the random
walk has the same asymptotic speed $\frac{\beta-1}{\beta+1}$ that we
would have seen on a naked $x$-axis without the traps. For $k \geq 0$
write $U(k)$
for the (random) time at which the random walk first arrives at the
vertex $(k,0)$. Establishing asymptotic speed $\frac{\beta-1}{\beta+1}$
is clearly the same as showing that
$\lim_{k \rightarrow \infty} \frac{U(k)}{k} = \frac{\beta+1}{\beta-1}$,
and for this, it is enough to show that
the time $\sum_{\{n: \, d_n < k\}} T_{n, tot}$ spent in traps
to the left of $(k, 0)$ satisfies
\begin{equation}  \label{eq:sum_of_delays_before_position_k}
\lim_{k \rightarrow \infty} \frac{1}{k} \sum_{\{n: \, d_n < k\}} T_{n, tot}
\, = \, 0 \quad \text{a.s.}
\end{equation}

By (\ref{eq:total_time_in_core}) and (\ref{eq:parameters_in_warmup_example}),
we have for all $n$ large enough that
\begin{eqnarray*}
\Pf[T_{n,tot}>T^*_{n,tot}]
& \leq & 2 \beta^{-e_n} \\
& \leq & 2 \beta^{- \alpha \log n}  \\
& = & 2 n^{- \alpha \log \beta} \, .
\end{eqnarray*}
Since $\beta>\beta_c = e^{1/\alpha}$ so that
$\alpha\log \beta>1$, we get
\[
\sum_{n=1}^\infty n^{- \alpha \log \beta} \, < \, \infty
\]
whence
\[
\sum_{n=1}^\infty \Pf[T_{n,tot}>T^*_{n,tot}]  \, < \, \infty \, .
\]
By Borel--Cantelli, we get a.s.\ that
$T_{n,tot}>T^*_{n,tot}$ for at most finitely many $n$, so that
\begin{equation}  \label{eq:by_Borel_Cantelli}
\sum_{n=1}^\infty (T_{n,tot}-T^*_{n,tot})  \, < \, \infty \,  \mbox{ a.s.}
\end{equation}
The left-hand side in (\ref{eq:sum_of_delays_before_position_k})
decomposes as
\begin{eqnarray}
\nonumber
\lim_{k \rightarrow \infty} \frac{1}{k} \sum_{\{n: \, d_n < k\}} T_{n, tot}
& = &
\lim_{k \rightarrow \infty} \frac{1}{k} \sum_{\{n: \, d_n < k\}}
(T_{n, tot} - T^*_{n,tot}) + T^*_{n,tot} \\
& = &
\label{eq:decomposed}
\lim_{k \rightarrow \infty} \frac{1}{k} \sum_{\{n: \, d_n < k\}}
(T_{n, tot} - T^*_{n,tot}) \\
& & +
\lim_{k \rightarrow \infty} \frac{1}{k} \sum_{\{n: \, d_n < k\}} T^*_{n,tot}
\nonumber
\end{eqnarray}
where the limit in (\ref{eq:decomposed}) is $0$ a.s.\ due to
(\ref{eq:by_Borel_Cantelli}). Hence, to settle the case $\beta> \beta_c$, it
suffices to show that
\[
\lim_{k \rightarrow \infty} \frac{1}{k} \sum_{\{n: \, d_n < k\}} T^*_{n,tot}
\, = \, 0
\]
or in other words that
\begin{equation}  \label{eq:remains_to_show_for_high_speed}
\lim_{n \rightarrow \infty} \frac{1}{d_n} \sum_{i=1}^n T^*_{i,tot}
\, = \, 0 \, .
\end{equation}
Lemma \ref{lem:exit_fast_if_core_not_hit} in combination with
(\ref{eq:geometric_distr}) yields
\[
\E[T^*_{i, tot}] \, \leq \, \frac{2\beta -1}{(\beta-1)^2} \, ,
\]
so that
\begin{eqnarray*}
\E \left[ \frac{1}{d_n} \sum_{i=1}^n T^*_{i,tot} \right]
& = & n^{-3} \sum_{i=1}^n \E[T^*_{i,tot}] \\
& = & \frac{2\beta -1}{n^2 (\beta-1)^2} \, .
\end{eqnarray*}
For any $\eps>0$, Markov's equality gives
\[
\Pf \left[ \frac{1}{d_n} \sum_{i=1}^n T^*_{i,tot} > \eps \right] \, \leq \,
\frac{2\beta -1}{\eps n^2 (\beta-1)^2} \, ,
\]
where we may note that the right-hand side is summable
over $n$, so that, by another application of Borel--Cantelli,
(\ref{eq:remains_to_show_for_high_speed}) follows,
and the proof of the proposition for $\beta>\beta_c$ is complete.

It remains to handle the case $\beta<\beta_c$. Write $A_n$ for
the event that the first time the random walk reaches $(d_n,0)$, it
immediately enters the trap and spends at least time $\beta^{c_n}=\beta^n$
in there. Note that $A_1, A_2, \ldots$ are independent, and
that, due to Lemma \ref{lem:stay_in_core},
\begin{eqnarray*}
\Pf[A_n] & \geq &
\frac{(\beta-1)^2}{(\beta+2)2\beta(\beta^{e_n+1}+\beta-2)} \\
& \geq & \frac{C}{\beta^{e_n}} \\
& \geq & \frac{C}{2\beta^{\alpha\log n}} \\
& = & \frac{C}{2n^{\alpha\log\beta}}
\end{eqnarray*}
for some $C>0$ which may depend on $\beta$ but not on $n$.

Next, for $i=1,2, \ldots$, define
\begin{equation}  \label{eq:def_Wi}
W_i \, = \, \sum_{n=2^{i-1}+1}^{2^i} I_{A_n}
\end{equation}
as the number of events happening amongst
$A_{2^{i-1}+1}, A_{2^{i-1}+2}, \ldots, A_{2^i}$. We get
\begin{eqnarray*}
\E[W_i] & = & \sum_{n=2^{i-1}+1}^{2^i} \Pf[A_n] \\
& \geq & 2^{i-1} \frac{C}{2 \cdot 2^{i \alpha \log \beta}} \\
& = & \frac{C}{4}2^{i(1- \alpha \log \beta)}
\end{eqnarray*}
where we may note that $1- \alpha \log \beta >0$ (this is where
we use $\beta<\beta_c$).
Note that for large $i$, $W_i$ is approximately Poisson, because
it counts independent events with small probabilities. Hence,
for $i$ large enough,
\begin{equation}  \label{eq:Poisson_Wi}
\Pf[W_i=0] \leq 2 \exp\left( - \frac{C}{4}2^{i(1- \alpha \log \beta)} \right)
\end{equation}
which decays to $0$ (faster than) exponentially, so that by Borel--Cantelli
we get a.s.\ that $W_i>0$ for all but at most finitely many $i$.
On the event that $W_i>0$, the time $U(2^{3i})$ of the first arrival of the
random walk at the vertex $(d_{2^i}, 0)=(2^{3i}, 0)$ satisfies
$U(2^{3i}) \geq \beta^{2^{i-1}}$. Hence we get a.s.\
that
\begin{equation}  \label{eq:very_late_arrivals}
\lim_{k \rightarrow \infty} \frac{U(k)}{k} = \infty
\end{equation}
along the subsequence $k=2^3, 2^6, 2^9, \ldots$. Since $U(k)$ is increasing,
$\frac{U(k)}{k}$ can drop by at most a factor $\frac{7}{8}$ as $k$
increases in the interval $[2^{3i},2^{3(i+1)})$, so convergence to $\infty$
along the full sequence in (\ref{eq:very_late_arrivals}) follows. This implies
zero asymptotic speed. $\Cox$

\section{The main construction}  \label{sect:main_construction}

In this section we specify the percolation process to be used as a witness
for proving Theorem \ref{th:main}. We proceed in three steps. First,
in Section \ref{sect:fractal}, we specify a (deterministic)
fractal-like percolation configuration that will play roughly
the same role as the path along the $x$-axis did in
Section \ref{sect:warm-up}.
Then, in Section \ref{sect:adding_traps}, we add traps to the construction.
Finally, in Section \ref{sect:translation_invariance}, we make the percolation
process translation invariant by means of a more successful application
of random translation than in Section \ref{sect:warm-up}.

\subsection{Fractal structure}  \label{sect:fractal}

For $k=1,2,\ldots$, a {\em branch} of order $k$ (see Figure 2) consists of
a horizontal path, called the {\em main part}, of length $b_k-1$, linked
at its rightmost vertex to a vertical path, called the {\em abutment},
of length $3\cdot2^{k-1}$, going
either upwards or downwards; here $(b_1, b_2, \ldots)$ is a fairly
rapidly growing sequence to be specified more precisely in
what follows. The point where the main part and the abutment meet is
called the {\em corner} of the branch, the other endpoint of the
main part is called the {\em tip}, and the other endpoint of the
abutment is called the {\em root}. The root of a branch of order $k$ will
always be situated somewhere on the main part of a branch of order $k+1$;
in this way, the random walk will be able to escape to infinity via
branches of higher and higher order.

\begin{figure}[h]
\centering{\mbox{\epsfig{file=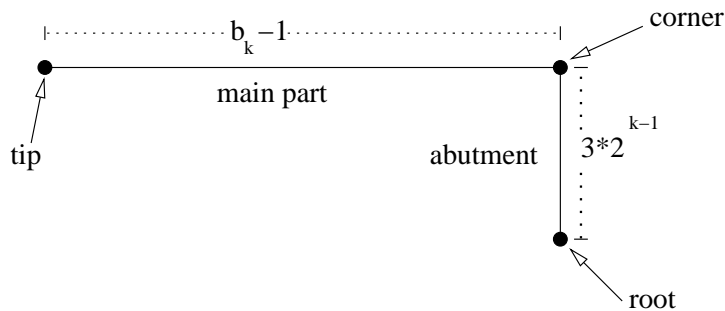,height=4cm}}}\caption{Branch of order $k$ with abutment pointing downwards.}
\end{figure}

Define a sequence $(q_1, q_2, \ldots)$ prescribing how many branches of order
$k$ should attach to each branch of order $k+1$. Each branch of
order $k+1$ will be attached to by $2q_k$ branches of order
$k$, $q_k$ of them attaching from above and the other $q_k$ from below.
In order for there to be room for all these branches of order $k$,
we will need the lengths to satisfy $b_{k+1} \geq q_k b_k$ for each $k$.
To give room for inserting traps in Section \ref{sect:adding_traps}, we will
need some extra margin, and will take
\begin{equation}  \label{eq:choice_of_b_k}
b_{k+1} \, = \, (q_k+1)b_k \, .
\end{equation}
A useful choice of $(q_1, q_2, \ldots)$ turns out to be
\[
q_k \, = \, 3^k -1
\]
so that bootstrapping (\ref{eq:choice_of_b_k}) gives
\begin{eqnarray*}
b_k & = & b_1 \prod_{i=1}^{k-1} (q_i+1)  \\
& = & b_1 \prod_{i=1}^{k-1} 3^i \\
& = & b_1 3^{k(k-1)/2} \, .
\end{eqnarray*}
Somewhat arbitrarily we set $b_1=4$, so that
\begin{equation}  \label{eq:final_choice}
b_k = 4 \cdot 3^{k(k-1)/2} \, .
\end{equation}

Here is how we arrange the branches. First, $y$-coordinates satisfying
$y=3l$ with integer $l$ will be reserved for (the main parts of)
branches. (Other $y$-coordinates will be used for traps later on.) More
specifically,
\[
\left\{
\begin{array}{c}
\mbox{$y$-coordinates with $y=3l$ with $l$ odd are reserved for branches
of order $1$} \\
\mbox{$y$-coordinates with $y=6l$ with $l$ odd are reserved for branches
of order $2$} \\
\vdots \\
\mbox{$y$-coordinates with $y=2^{k-1}3l$ with $l$ odd are reserved for branches
of order $k$} \\
\vdots
\end{array} \right.
\]

\begin{figure}[h]
\mbox{\epsfig{file=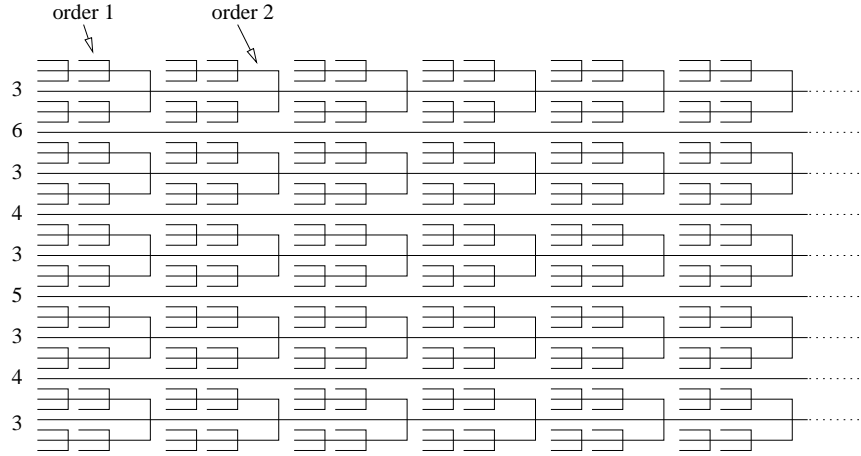,height=6cm}}\caption{Self-similar structure of branches. The orders of the larger branches are indicated to the left.}
\end{figure}

\noindent Now, for each $k$ and each $y$ with $y=2^{k-1}3l$ for $l$ odd, we set up a
branch of order $k$ with its tip at $(0, y)$, its corner at $(b_k-1, y)$,
and its root at either $(b_k-1, y + 3 \cdot 2^{k-1})$ (i.e., the abutment
pointing upwards)) or
$(b_k-1, y - 3 \cdot 2^{k-1})$ (the abutment pointing downwards)
chosen as follows. Since $y=2^{k-1}3l$ with $l$ odd, we have that
exactly one of the numbers $y + 3 \cdot 2^{k-1}$ and $y - 3 \cdot 2^{k-1}$
equals $2^k3j$ for some $j$ odd, and we choose to direct the abutment so that
the root ends up at such a $y$-coordinate, thus ensuring that it
sits on (the main part of) a branch of order $k+1$.

In this way, for any $k\geq 2$, any branch of order $k$ having its tip
at, say, $(x,y)$ will have
received exactly two branches of order $k-1$ attaching to it,
one from above and one from below, and both of them with the root at
$(x+b_{k-1}-1, y)$. For each such branch of order $k$, we attach another
$2(q_{k-1}-1)$ branches of order $k-1$ to it, two
(one from above and one from below) at each of the points
$(x+2b_{k-1}-1, y), (x+3b_{k-1}-1, y), \ldots, (x+q_{k-1}b_{k-1}-1, y)$.

This creates, for any $k\geq 2$, further branches of order $k$ to which
presently there are no branches of order $k-1$ attached. To such
a branch of order $k$ with tip at $(x,y)$
we attach $2q_{k-1}$ branches of order $k-1$, two
(one from above and one from below) at each of the points
$(x+b_{k-1}-1, y), (x+2b_{k-1}-1, y), \ldots, (x+q_{k-1}b_{k-1}-1, y)$.

Repeating this procedure {\em ad infinitum} produces a percolation
configuration $\xi \in \{{\rm open,closed}\}^E$, where as before
$E$ is the edge set
of the square lattice. This is the fractal structure that forms the
foundation of our construction.

\subsection{Adding traps}   \label{sect:adding_traps}

Now we will add traps to the configuration $\xi$ of
Section \ref{sect:fractal}. Each branch of order $k \geq 2$ will be
equipped with exactly one trap. The trap will be situated in the
final region of the main part of the branch. This part contains a stretch
of length $(b_k-1)- (q_{k-1}b_{k-1}-1)$ from the last attachment of a
lower-order branch to its corner point. By (\ref{eq:choice_of_b_k})
and (\ref{eq:final_choice}),
this length equals
\begin{eqnarray*}
(b_k-1)- (q_{k-1}b_{k-1}-1) & = & b_{k-1} \\
& = & 4 \cdot 3^{(k-1)(k-2)/2}  \, .
\end{eqnarray*}
We choose to place the anchor of the trap exactly at the midpoint
of this stretch, i.e., exactly $2 \cdot 3^{(k-1)(k-2)/2}$ steps to
the left of the corner point of the order $k$-branch. The
lengths $e_k$ and $c_k$ of the trap's entrance and core, respectively,
are preliminarily chosen as
\begin{equation}  \label{eq:preliminary_choice}
e_k \, = \, \left\lceil \frac{\log k}{\log \gamma} \right\rceil
\end{equation}
and
\[
c_k \, = \, 3^{(k-1)(k-2)/2}
\]
respectively. For values of $\gamma$ very close to $1$, it may turn out that
the chosen value of $e_k$ in (\ref{eq:preliminary_choice}) exceeds
$2 \cdot 3^{(k-1)(k-2)/2}$ so that the trap bumps into
the abutment of the last branch of order $k-1$ attaching from above to the
branch of order $k$ to which the trap is attached. This cannot be allowed
to happen, and we therefore replace (\ref{eq:preliminary_choice}) by
\begin{equation}  \label{eq:corrected_choice}
e_k \, = \, \min \left\{ \left\lceil \frac{\log k}{\log \gamma} \right\rceil,
3^{(k-1)(k-2)/2} \right\} \, ,
\end{equation}
and note that this coincides with (\ref{eq:preliminary_choice}) for
all $k$ large enough.

Write $\eta \in \{{\rm open,closed}\}^E$ for the percolation
configuration obtained by adding traps in this manner
to the configuration $\xi$.

\subsection{Stationarizing}  \label{sect:translation_invariance}

Let $\Psi''$ be the probability measure on $\{{\rm open,closed}\}^E$
corresponding to picking the configuration $\eta \in \{{\rm open,closed}\}^E$
from Section \ref{sect:adding_traps} deterministically.
$\Psi''$ is turned into a vertically translation invariant
measure $\Psi'$ on $\{{\rm open,closed}\}^E$ by shifting the configuration
$\eta$ vertically by an amount $S_y$ chosen uniformly from
$\{-m, -m+1, \ldots, m\}$, and taking weak limits as $m \rightarrow \infty$,
if necessary after passing to a subsequence. It is readily checked that,
for any $k \geq 1$,
\begin{equation}  \label{eq:vertical_nondisappearance}
\Psi'(\mbox{the origin is the tip of a branch of order $k$})
\, = \, \frac{1}{3 \cdot 2^k} \, .
\end{equation}

To achieve translation invariance also in the horizontal direction, we
first choose a configuration $\eta'\in \{{\rm open,closed}\}^E$
according to $\Psi'$,
and then we pick a shift $S_x \in \{0, \ldots, b_n-1\}$ independently of $\eta'$, resulting in a
probability measure $\Psi'_n$ on $\{{\rm open,closed}\}^E$.
We finally take the probability measure $\Psi$ on $\{{\rm open,closed}\}^E$
as a weak limit of the $\Psi'_n$ measures as $n \rightarrow \infty$,
if necessary after passing to a subsequence.\footnote{In fact, it
turns out that the limits
exist without passing to a subsequence, both here and in
going from $\Psi''$ to $\Psi'$, but we do not need this.}

It is clear that $\Psi$ is both vertically and horizontally
translation invariant.
We now need to check that local structures (branches
of order $k$ for given $k$, and traps) do not disappear upon us in
the limit as in
the failed attempt at stationarizing in Section \ref{sect:warm-up}.
For this, it suffices to show that for any $k \geq 1$,
\begin{equation}  \label{eq:need_to_show}
\lim_{n \rightarrow \infty} \Psi'_n (\mbox{the origin sits on the
main part of a branch of order $k$}) \, > \, 0 \, .
\end{equation}
From (\ref{eq:vertical_nondisappearance}), we get immediately that
\begin{equation}  \label{eq:y-coordinate}
\Psi'_n (\mbox{the origin is at a $y$-coordinate devoted
to branches of order $k$})
\, = \, \frac{1}{3 \cdot 2^k} \, .
\end{equation}
Writing $A_k$ for the event in (\ref{eq:y-coordinate}), we get
for $n>k$ using a direct count of the $b_n$ different horizontal
translations available to $\Psi'_n$ that
\begin{eqnarray*}
\lefteqn{\hspace{-25mm}
\Psi'_n(\mbox{the origin is on the main part of a branch of order $k$}
\, | \, A_k)} \\
& = & \frac{b_k \prod_{i=k}^{n-1} q_i}{b_n} \, = \,
\frac{b_k \prod_{i=k}^{n-1} q_i}{b_k \prod_{i=k}^{n-1} (q_i+1)} \\
& = & \prod_{i=k}^{n-1} \frac{q_i}{q_i+1}
\, = \, \prod_{i=k}^{n-1} (1-3^{-i}) \, ,
\end{eqnarray*}
where the second equality derives from the
recursive definition (\ref{eq:choice_of_b_k}).
Multiplying by (\ref{eq:y-coordinate}) and sending $n \rightarrow \infty$
gives
\begin{eqnarray*}
\lefteqn{\hspace{-32mm}
\Psi(\mbox{the origin is on the main part of a branch of order $k$})} \\
& = & \frac{1}{3 \cdot 2^k}
\prod_{i=k}^\infty(1-3^{-i}) \, > \, 0 \, ,
\end{eqnarray*}
so that the local structures do not disappear upon us in the limit.

\section{Positive speed without the traps}  \label{sect:without_traps}

In this section, we study what happens to the random walk in the modified
percolation process obtained by removing all the traps. To this end,
write $\Psi^*$ for the probability measure on $\{{\rm open,closed}\}^E$
corresponding to picking a percolation configuration according to $\Psi$
and then deleting all the traps. Similarly as in Theorem \ref{th:main},
we write $\Pb_{\Psi^*,\beta}$ for the joint law of the percolation
configuration chosen according to $\Psi^*$ and the random
walk $\{Z_t\}_{t \geq 0}$ with drift
parameter $\beta$ starting at the origin.
Recall that $\{0 \leftrightarrow \infty\}$ denotes the
event that the origin is
in the infinite cluster of
the percolation configuration.
\begin{prop}   \label{prop:without_traps}
For any $\beta > 1$ we have
\[
\lim_{t \rightarrow \infty} \frac{X_t}{t} \, = \, \frac{\beta-1}{\beta+1}
\, \, \, \Pb_{\Psi^*,\beta} \mbox{-a.s. on the event }
\{0 \leftrightarrow \infty\} \, .
\]
\end{prop}
A first simplification for the proof of Proposition \ref{prop:without_traps}
is the following reduction, where we write $B_k$ for the event that
the origin is on the main part of a branch of order $k$.
\begin{lemma}  \label{lem:reduction}
Suppose for given $\beta>1$ and $\theta\geq 0$
that
\begin{equation}  \label{eq:from_the_special_case...}
\lim_{t \rightarrow \infty} \frac{X_t}{t} \, = \, \theta \, \, \,
\Pb_{\Psi^*,\beta} \mbox{-a.s. on the event } B_1 \, .
\end{equation}
Then, in fact,
\begin{equation}  \label{eq:...to_the_general}
\lim_{t \rightarrow \infty} \frac{X_t}{t} \, = \, \theta \, \, \,
\Pb_{\Psi^*,\beta} \mbox{-a.s. on the event } \{0 \leftrightarrow \infty\}.
\end{equation}
The same result holds with $\Pb_{\Psi,\beta}$
in place of  $\Pb_{\Psi^*,\beta}$.
\end{lemma}
{\bf Proof.} Fix a probability distribution ${\bf Q}$ on $\Z^2$ with
full support, and pick $\tilde{z}=(\tilde{x}, \tilde{y}) \in \Z^2$
according to ${\bf Q}$. Given the percolation configuration chosen according
to $\Psi^*$, run {\em two}
random walks $\{Z_t\}_{t \geq 0} = \{(X_t, Y_y)\}_{t \geq 0}$ and
$\{\tilde{Z}_t\}_{t \geq 0}=
\{(\tilde{X}_t, \tilde{Y}_t)\}_{t \geq 0}$
starting at the origin and at $\tilde{z}$, respectively,
coupled as follows. Let $\{\tilde{Z}_t\}_{t \geq 0}$ run independently
of $\{Z_t\}_{t \geq 0}$ except that from the first time $T$ that
$\tilde{Z}_t$ hits the origin (if ever), $\tilde{Z}_t$ plagiarizes the
trajectory of $Z_t$ from then on, meaning that
$\{\tilde{Z}_T, \tilde{Z}_{T+1}, \ldots\} = \{Z_0, Z_1, \ldots\}$.
By translation invariance of $\Psi^*$,
\begin{equation}  \label{eq:translation_invariance}
\{\tilde{Z}_t - \tilde{z} \}_{t \geq 0} \, \, \, \mbox{has the same
distribution as} \, \, \, \{Z_t\}_{t \geq 0} \, .
\end{equation}
Define $D= D_1 \cap D_2 \cap D_3$ where
\[
D_1 \, = \, \{ 0 \leftrightarrow \infty \} \cap
\left\{\lim_{t \rightarrow\infty} \frac{X_t}{t} \neq \theta \right\} \, ,
\]
\[
D_2 \, = \, \{\tilde{z} \mbox{ is on the main part of an order-$1$ branch}\}
\, ,
\]
and $D_3 = \{ T < \infty\}$. Assume for contradiction that
(\ref{eq:from_the_special_case...}) holds and that
(\ref{eq:...to_the_general}) fails. Then
$D_1$ has positive probability, and it is easy too see that in that case,
$D= D_1 \cap D_2 \cap D_3$ has positive probability too. On the
event $D$ we get that
\[
\lim_{t \rightarrow \infty}\frac{\tilde{X}_t - \tilde{x}}{t} \, = \,
\lim_{t \rightarrow \infty}\frac{X_{t-T} - \tilde{x}}{t} \, = \,
\lim_{t \rightarrow \infty}\frac{X_t}{t} \, \neq \, \theta
\]
which, in view of (\ref{eq:translation_invariance}), contradicts
(\ref{eq:from_the_special_case...}). This concludes the argument for
$\Pb_{\Psi^*,\beta}$, and the same argument goes through with
$\Pb_{\Psi,\beta}$
in place of  $\Pb_{\Psi^*,\beta}$. $\Cox$

\medskip\noindent
For the purpose of proving Proposition \ref{prop:without_traps}, we may
(due to Lemma \ref{lem:reduction}) assume that the origin sits on
the main part of a branch of order $1$, and go on to analyze random
walk from there. In this case, there is a unique self-avoiding
path $P$ in the percolation configuration from the origin to infinity.
This path goes through (parts of) the main parts and the abutments of
branches of increasing order $1,2,3, \ldots$. Write $a_k$ for the
$x$-coordinate of the abutment of the branch of order $k$ in this path.
For $k \geq 2$, a crude lower bound
for $a_k$, which follows directly from the construction, is
\begin{equation}   \label{eq:lower_bound_on_a_k}
a_k \, \geq \, b_{k-1} \, = \, 4 \cdot 3^{(k-1)(k-2)/2}  \, .
\end{equation}

A couple of further lemmas will be convenient to isolate for the
proof of  Proposition \ref{prop:without_traps}.
\begin{lemma}  \label{lem:escape_probability}
For any $k$ and any vertex $z$ which sits on the main part of a branch
of order
$k$ and which sits at least $3 \cdot 2^{k-1}$ steps to the left of the corner
point of that branch, we have that the escape probability for the random walk
starting from $z$ is at least
\[
\frac{\beta-1}{2(3+ \beta)} \, .
\]
\end{lemma}

\medskip\noindent
{\bf Proof.} We proceed electrically as in Section \ref{sect:warm-up},
attaching a resistance $R(e)= \beta^{-x(e)}$ to each edge $e$ of
the percolation configuration, where $x(e)$ is the largest
$x$-coordinate amongst the two endpoints of $e$.
We recall from (\ref{eq:escape_probability}) that the escape
probability from a vertex $z=(x,y)$
is given by
\begin{equation}  \label{eq:escape_probability_once_more}
p_{esc}(z) \, = \, \left( R_{eff}(z, \infty)
\sum_{\{e: \, e \sim z\}} R(e)^{-1}   \right)^{-1}
\end{equation}
where the sum is over all edges $e$ that are incident to vertex $z$.
The sum in (\ref{eq:escape_probability_once_more}) is bounded above by
the sum
\begin{equation}   \label{eq:conductances_in_neighborhood}
3\beta^x+ \beta^{x+1} \, ,
\end{equation}
corresponding to the case where all four possible edges incident to $z$ are
present in the percolation configuration.
The effective resistance from $z$ to $\infty$ is simply the sum of
the resistances along the unique self-avoiding path from
$z$ to $\infty$. Counting only the horizontal edges of this path would give simply the sum $\sum_{i=x+1}^\infty \beta^{-i}$. The point of the choice
of the bound $3 \cdot 2^{k-1}$ in the lemma is that the set of vertical
edges on the path can be paired with a subset of the set of horizontal edges
on the path, in such a way that a vertical edge is always paired with
a horizontal edge with smaller $x$-coordinate and therefore
larger resistance. Hence the set of vertical edges can contribute at most
as much as the set of horizontal edges to $R_{eff}(z, \infty)$,
so
\[
R_{eff}(z, \infty) \,  \leq 2 \, \sum _{i=x+1}^\infty \beta^{-i}
\, = \, \frac{2\beta^{-x}}{\beta-1} \, .
\]
Plugging this bound and (\ref{eq:conductances_in_neighborhood}) into
(\ref{eq:escape_probability_once_more}) gives
\begin{eqnarray*}
p_{esc}(z) & \geq & \left( (3\beta^x+ \beta^{(x+1)})
\frac{2\beta^{-x}}{\beta-1} \right)^{-1} \\
& = &  \frac{\beta-1}{2(3+ \beta)}
\end{eqnarray*}
as desired. $\Cox$

\begin{lemma}  \label{lem:expected_duration_of_excursion}
For any $\beta>1$, there exists a constant $C'= C'_\beta$ independent of
$k$, such that
\begin{description}
\item{(a)} a random walk taking a step to the left from a
corner point of a branch of order $k$ has an expected time
until return to the corner point which
is at most $C'_\beta$, and
\item{(b)} a random walk taking a step into a branch of order $k$ from its
root has an expected return time to the root which is at most
\[
3 \cdot 2^k + C'_\beta \, .
\]
\end{description}
\end{lemma}

\medskip\noindent
{\bf Proof.}
Imagine random walk with bias $\beta$ on a finite connected subgraph
$\tilde{G}=(\tilde{V}, \tilde{E})$
of the
square lattice. This can be described as a finite-state Markov chain with
a unique stationary distribution, where it is easily checked that each vertex
$v \in \tilde{V}$ receives a probability $\pi(v)$ proportional to
the sum
\begin{equation}  \label{eq:conductance_of_vertex}
\sum_{e \sim v} R(e)^{-1}
\end{equation}
of inverse edge resistances (defined in
the same way as in the proof of Lemma \ref{lem:escape_probability}) among
edges incident to $v$.
Define $R^{-1}(v)$ as the sum in (\ref{eq:conductance_of_vertex}).
Standing at a given $v \in \tilde{V}$, the expected return time to $v$
is
\begin{equation}  \label{eq:expected_return_time_general_formula}
\pi(v)^{-1} \, = \, \frac{1}{R^{-1}(v)}\sum_{w \in \tilde{V}}R^{-1}(w) \, .
\end{equation}
For a vertex $z =(x,y)$, we have that
\begin{equation}  \label{eq:inverse_resistance_of_a_node}
R^{-1}(z) \, \leq \, 3\beta^x+ \beta^{x+1} \, .
\end{equation}
For our percolation process, formula
(\ref{eq:expected_return_time_general_formula})
applies when the random walk leaves
a vertex $v$ to enter a finite region of the percolation configuration
cut of from $v$ from the rest of the configuration.  Applying this when
$z=(x,y) $ is a corner point of a
branch of order $k$, we have that the entire
finite structure $\tilde{G}=(\tilde{V}, \tilde{E})$
cut off by $z$ is contained in the cone
\begin{equation}  \label{eq:cone}
\{z'= (x', y') \in \Z^2: \, x' \leq x, y' \in [y-(x-x'), y+ (x-x')] \} \, .
\end{equation}
Summing (\ref{eq:inverse_resistance_of_a_node}) over this cone gives
\[
\sum_{w \in \tilde{V}} R^{-1}(w) \, \leq \,
\sum_{i=0}^\infty (2i+1)(3\beta^{x-i}+ \beta^{x-i+1}) \, .
\]
Furthermore, the $R^{-1}$ value of the corner point itself is
$R^{-1}(z)= 2 \beta^x$. Plugging these observations into
(\ref{eq:expected_return_time_general_formula}) yields that the
expected return time to $z$ is bounded by
\begin{equation}  \label{eq:sum_in_cone}
\frac{1}{2 \beta^x} \sum_{i=0}^\infty (2i+1)(3\beta^{x-i}+ \beta^{x-i+1})
\, = \, \frac{3 + \beta}{2} \sum_{i=0}^\infty (2i+1)\beta^{-i}\,  < \infty \, ,
\end{equation}
so part (a) of the lemma is established with $C'_\beta$ equal to the
right-hand side in (\ref{eq:sum_in_cone}).

Part (b) follows similarly, the only difference being that
in the sum in (\ref{eq:expected_return_time_general_formula})
we have to
take into account the additional $3 \cdot 2^{k-1}$ vertices
in the abutment, each of which contributes an amount $2\beta^x$ to
the sum.
$\Cox$

\medskip\noindent
{\bf Proof of Proposition \ref{prop:without_traps}.}
To establish that $\lim_{t \rightarrow \infty} \frac{X_t}{t} =
\frac{\beta-1}{\beta+1}$, it suffices to show that
\begin{equation}  \label{eq:inverse_asymptotic_speed}
\lim_{x \rightarrow \infty} \frac{U(x)}{x} \, = \, \frac{\beta+1}{\beta - 1}
\end{equation}
where, similarly as in the proof of
Proposition \ref{prop:warm-up_construction}, we define $U(x)$ as the time
of first arrival at $x$-coordinate $x$:
\[
U(x) \, = \, \min \{t: \, X_t=x \} \, .
\]
As a means towards estimating $U(x)$ well enough to establish
(\ref{eq:inverse_asymptotic_speed}), we decompose it as
\[
U(x) \, = \, U'(x) + U''(x) \,
\]
where $U'(x)$ is the time spent on the path $P$ to infinity before hitting $x$-coordinate
$x$, and $U''(x)$ is the time spent outside the path $P$ before first
hitting $x$-coordinate $x$. Our plan is to establish
\begin{equation}  \label{eq:U-prime}
\lim_{x\rightarrow\infty} \frac{U'(x)}{x} \, = \, \frac{\beta+1}{\beta-1}\quad\text{a.s.}
\end{equation}
and
\begin{equation}  \label{eq:U-biss}
\lim_{x\rightarrow\infty} \frac{U''(x)}{x} \, = \, 0 \quad\text{a.s.}
\end{equation}
Together, (\ref{eq:U-prime}) and (\ref{eq:U-biss}) will of course imply
(\ref{eq:inverse_asymptotic_speed}).

We begin with (\ref{eq:U-prime}). Note that for the purpose of
studying $U'(x)$ we may simply pretend that paths other than $P$ do
not exist. Now, if $P$ had no abutments and instead consisted of a single
straight path along the $x$-axis, (\ref{eq:U-prime}) would follow
immediately from the strong law of large numbers. Furthermore, it
is easy to see that $\{U'(x)\}_{x \geq 1}$ stochastically dominates
the corresponding process in such an ideal scenario. Hence
\begin{equation}  \label{eq:half-way-there}
\lim_{x \rightarrow \infty} \frac{U'(x)}{x} \, \geq \,
\frac{\beta+1}{\beta-1}\, .
\end{equation}
To strengthen this to an equality, we need to show that the delay caused
by abutments is small. More precisely, define $S(k)$ as the time spent
on the abutment of order $k$ in $P$,
{\em plus} the time spent on $P$ to the left of this abutment
after first having visited its corner point. The inequality
(\ref{eq:half-way-there}) is strengthened to the equality
(\ref{eq:U-prime}) if we can show that
\[
\lim_{x \rightarrow \infty} \frac{1}{x} \sum_{\{k: a_k \leq x\}} S(k)
\, = \, 0 \quad\text{a.s.}
\]
which is equivalent to
\begin{equation}  \label{eq:this_subsequence_is_enough}
\lim_{k \rightarrow \infty} \frac{1}{a_k} \sum_{j=1}^k S(j) \, = \, 0\quad\textrm{a.s.}
\end{equation}
We go on to estimate $\E[S(k)]$. On the abutment itself, the
random walk behaves like simple
random walk, and it is a standard fact that the expected time
on it until first hitting the
root point equals its length squared, i.e., $(3 \cdot 2^{k-1})^2$.
During this walk, the expected number of
returns to the corner point is linear in the length $3 \cdot 2^{k-1}$, and
to each
such return corresponds a geometric number with mean $1$
of excursions
to the left of the corner point, so that the expected total number of
such excursions is again linear in the abutment's length. Lemma
\ref{lem:expected_duration_of_excursion} (a) ensures that the expected
duration of such an excursion is bounded uniformly in $k$.
Furthermore, Lemma \ref{lem:escape_probability} ensures that
once the walk has reached the root point,
the number of times it goes back into the abutment again is dominated by a geometric variable with mean $(\frac{\beta-1}{2(3+\beta)})^{-1}$, and Lemma
\ref{lem:expected_duration_of_excursion} (b) ensures
that each such excursion has expected duration at most
$3 \cdot 2^k +C'_\beta$.
Summing up the contributions to $S(k)$, we get that
there exists a constant $C=C_\beta$ independent of $k$ such that
\[
\E[S(k)] \, \leq \, C_\beta 2^{2k} \, .
\]
Taking expectation in the left-hand side of
(\ref{eq:this_subsequence_is_enough}) and plugging in
(\ref{eq:lower_bound_on_a_k}) gives
\begin{eqnarray*}
\E\left[\frac{1}{a_k} \sum_{j=1}^k S(j) \right]
& \leq & \frac{1}{4 \cdot 3^{(k-1)(k-2)/2}} \sum_{j=1}^k \E[S(j)] \\
& \leq & \frac{C_\beta}{4 \cdot 3^{(k-1)(k-2)/2}} \sum_{j=1}^k 2^{2j} \\
& \leq & \frac{C_\beta k 2^{2k}}{4 \cdot 3^{(k-1)(k-2)/2}} \, .
\end{eqnarray*}
Markov's inequality gives, for any $\eps>0$, that
\[
\Pf \left[ \frac{1}{a_k} \sum_{j=1}^k S(j) \geq \eps \right] \, \leq \,
\frac{C_\beta k 2^{2k}}{4 \cdot 3^{(k-1)(k-2)/2}\eps}
\]
which decays to $0$ exponentially fast as $k \rightarrow \infty$,
so Borel--Cantelli gives (\ref{eq:this_subsequence_is_enough}). Hence,
(\ref{eq:U-prime}) is established, and it only remains to prove
(\ref{eq:U-biss}).

For $k \geq 2$, define $W_k$ as the total time spent in parts of
the percolation configuration away from $P$ that attach to $P$ in the
part of $P$ that belongs to a branch of order $k$. Regions contributing
to $W_k$ are of two kinds, namely,
\begin{description}
\item{(i)} branches of order $k-1$ (together with their
respective subbranches), and
\item{(ii)} the section not contained in $P$ of the branch of order $k$
itself.
\end{description}
There are at most $2q_k= 2(3^k-1)$ branches of order $k-1$
contributing to $W_k$. By Lemma \ref{lem:escape_probability}, the
expected number of times that each such branch
is visited is at most $(\frac{\beta-1}{2(3+\beta)})^{-1}$, and by Lemma
\ref{lem:expected_duration_of_excursion} (b) the expected
duration of each such visit is at most $3 \cdot 2^k + C'_\beta$. Hence,
the contribution from (i) to $\E[W_k]$ is at most
\[
2(3^k-1)\left(\frac{\beta-1}{2(3+\beta)}\right)^{-1}(3 \cdot 2^k + C'_\beta) \, .
\]
The contribution from (ii) is obtained by multiplying the expected number
of visits to the section in question, by the expected
duration of each visit; the former is bounded by
$(\frac{\beta-1}{2(3+\beta)})^{-1}$ due to Lemma \ref{lem:escape_probability},
and the latter is bounded by $C'_\beta$ by arguing as in the proof of Lemma
\ref{lem:expected_duration_of_excursion} (a). Summing the contributions
from (i) and (ii) gives
\[
\E[W_k] \, \leq \, \left(\frac{\beta-1}{2(3+\beta)}\right)^{-1}
\left( 2(3^k - 1) (3 \cdot 2^k + C'_\beta) + C'_\beta \right) \, .
\]
For $k \geq 3$ and with $a_k$ as in (\ref{eq:lower_bound_on_a_k}), we get
\begin{eqnarray*}
\lefteqn{
\E\left[ \frac{W_2 + W_3 + \cdots + W_k}{a_{k-1}} \right]} \\
& \leq &
\frac{1}{4 \cdot3^{(k-2)(k-3)/2}} \sum_{j=2}^k
\left(\frac{\beta-1}{2(3+\beta)}\right)^{-1}
\left( 2(3^j - 1) (3 \cdot 2^j + C'_\beta) + C'_\beta \right) \\
& \leq &
\frac{k-1}{4 \cdot3^{(k-2)(k-3)/2}}
\left(\frac{\beta-1}{2(3+\beta)}\right)^{-1}
\left( 2(3^k - 1) (3 \cdot 2^k + C'_\beta) + C'_\beta \right)
\end{eqnarray*}
which tends to $0$ exponentially fast in $k$. Hence, by Markov's inequality
and Borel--Cantelli, a.s.\ $(W_2 + \cdots + W_k)/a_{k-1}$ will exceed
any given $\eps>0$ at most finitely many times. In other words,
we have a.s.\ that
\begin{equation}  \label{eq:almost_done}
\lim_{k \rightarrow \infty} \frac{W_2 + W_3 + \cdots + W_k}{a_{k-1}}
\, = \, 0 \, .
\end{equation}
Now, it is easy to see that $\frac{U''(x)}{x} \leq
\frac{W_2 + W_3 + \cdots + W_{k'}}{a_{k'-1}}$ where $k'$ is the smallest
$k$ such that $a_k \geq x$. Hence (\ref{eq:almost_done}) implies the
desired (\ref{eq:U-biss}), so the proof is complete. $\Cox$

\section{Main construction: positive speed regime}
\label{sect:positive_speed}

We are almost ready to switch from considering the modified percolation
configuration gotten from $\Psi^*$, to the full percolation configuration,
including traps, obtained from $\Psi$. But before taking the
full step we make an intermediate stop at the probability measure
$\Psi^{**}$ on $\{{\rm open,closed}\}^E$ corresponding to
picking a configuration according to $\Psi$ and then deleting all
traps situated directly on the path $P$ from $0$ to $\infty$, but leaving
all other traps undeleted. We have the following variation of
Proposition \ref{prop:without_traps}.
\begin{prop}   \label{prop:without_some_of_the_traps}
For any $\beta > 1$ we have
\[
\lim_{t \rightarrow \infty} \frac{X_t}{t} \, = \, \frac{\beta-1}{\beta+1}
\, \, \, \Pb_{\Psi^{**},\beta} \mbox{-a.s. on the event }
\{0 \leftrightarrow \infty\} \, .
\]
\end{prop}
{\bf Proof.}
The proof of Proposition \ref{prop:without_traps} translates verbatim
to this case. The crucial point to note is that the estimates in
Lemma \ref{lem:expected_duration_of_excursion} for the time
spent in branches outside of $P$ are still valid when
traps are added, because any trap added to such a branch
(or any of its subbranches) will be contained in the cone
(\ref{eq:cone}). $\Cox$

\medskip\noindent
In this section we consider the large drift regime $\beta>\beta_c$. In
view of Proposition \ref{prop:without_some_of_the_traps}, all we need
to keep track of is the time spent in the traps directly attached to the
path $P$. The trap attached to the order-$k$ branch part of
$P$ will henceforth be called {\em trap number $k$}.
We go on to consider random walk on the full percolation
configuration obtained from $\Psi$.
Define $U'''(x)$ as the time spent in traps directly attached
to the path $P$ before first hitting $x$-coordinate $x$. By
reasoning similarly
as in the decomposition of $U(x)$ at the beginning of the proof
of Proposition \ref{prop:without_traps}, what we
need to show is that a.s.\
\begin{equation}  \label{eq:triss}
\lim_{x \rightarrow \infty} \frac{U'''(x)}{x} \, = \, 0 \, .
\end{equation}
Analogously to the notation in Section
\ref{sect:warm-up}, we write (for $k \geq 2$ and $i \geq 1$)
$T_{k,i}$ for the time
spent in the trap attached to the order-$k$ branch
(trap number $k$, for short)
in $P$ during
the $i$'th visit to this trap; if $i$ exceeds the number of visits
to the trap, we set $T_{k,i}=0$. We also define the total
time spent in the trap
\[
T_{k, tot} \, = \, \sum_{i=1}^\infty T_{k,i} \, .
\]
Still following Section \ref{sect:warm-up}, we define
\[
T^*_{k,i} \, = \left\{
\begin{array}{ll}
0 & \mbox{if the walk hits the trap's core during this visit} \\
T_{k,i} & \mbox{otherwise,}
\end{array} \right.
\]
and $T^*_{k, tot} = \sum_{i=1}^\infty T^*_{k,i}$.

\medskip\noindent
{\bf Proof of Theorem \ref{th:main}, case $\beta>\beta_c$.}
By Lemma \ref{lem:reduction}, we may assume that the origin
sits on a branch of order $1$.
We begin by noting that
\[
\frac{U'''(x)}{x} \, \leq \, \frac{1}{a_{k'-1}}\sum_{j=2}^{k'} T_{j,tot}
\]
where $k'$ is the smallest
$k$ such that $a_k \geq x$. Hence, to establish the desired
(\ref{eq:triss}), it suffices to show
that a.s.
\begin{equation}  \label{eq:need_to_show_for_large_bias}
\lim_{k \rightarrow\infty} \frac{1}{a_{k-1}}\sum_{j=2}^{k} T_{j,tot}
\, = \, 0 \, .
\end{equation}
Next, we note that the expected number of times
that trap number $k$ is visited is at most
$(\frac{\beta-1}{2(3+\beta)})^{-1}$ due to
Lemma \ref{lem:escape_probability}. In combination with
Lemma \ref{lem:exit_fast_if_core_not_hit}, this gives
\[
\E[T^*_{j, tot}] \, \leq \,
\left(\frac{\beta-1}{2(3+\beta)}\right)^{-1}
\left(\frac{2\beta-1}{\beta-1}\right)
\]
and, using (\ref{eq:lower_bound_on_a_k}),
\begin{eqnarray*}
\E \left[
\frac{1}{a_{k-1}}\sum_{j=2}^{k} T^*_{j,tot} \right] & \leq &
\frac{k \left(\frac{\beta-1}{2(3+\beta)}\right)^{-1}
\left(\frac{2\beta-1}{\beta-1}\right)}{ 4 \cdot 3^{(k-1)(k-2)/2}}
\end{eqnarray*}
which decays to $0$ (faster than) exponentially as $k \rightarrow \infty$.
This allows us to exploit the familiar combination of Markov's inequality
and Borel--Cantelli: for any $\eps>0$ the probability that
$\frac{1}{a_{k-1}}\sum_{j=2}^{k} T^*_{j,tot}$ exceeds $\eps$ is summable
over $k$, so that a.s.
\begin{equation}  \label{eq:need_to_show^*}
\lim_{k \rightarrow\infty} \frac{1}{a_{k-1}}\sum_{j=2}^{k} T^*_{j,tot}
\, = \, 0 \, .
\end{equation}
This will imply the desired (\ref{eq:need_to_show_for_large_bias}) as soon
as we can establish that $T_{j, tot} > T^*_{j, tot}$ for at most
finitely many $j$. For this we proceed as in
the proof of Proposition \ref{prop:warm-up_construction}: Lemma
\ref{lem:hit_the_core} tells us that each time the
walk enters trap number $j$, it has probability
$\frac{\beta-1}{\beta^{e_j+1} + \beta -2}$ of hitting the core. Using
again that the expected number of visits to the trap is at most
$(\frac{\beta-1}{2(3+\beta)})^{-1}$, we get that
\begin{equation}  \label{eq:BC_estimate}
\Pf[T_{j, tot} > T^*_{j, tot}] \, \leq \,
\left(\frac{\beta-1}{2(3+\beta)}\right)^{-1}
\left(\frac{\beta-1}{\beta^{e_j+1} + \beta -2} \right) \, .
\end{equation}
The choice (\ref{eq:corrected_choice}) of $e_j$ gives
$e_j \geq \frac{\log j}{\log \beta_c}$, so that the
estimate (\ref{eq:BC_estimate}) may be further bounded as
\begin{eqnarray*}
\Pf[T_{j, tot} > T^*_{j, tot}] & \leq &
\left(\frac{\beta-1}{2(3+\beta)}\right)^{-1}
\left(\frac{\beta-1}{\beta^{e_j+1} + \beta -2} \right) \\
& \leq &
\left(\frac{\beta-1}{2(3+\beta)}\right)^{-1}
\left(\frac{\beta-1}{\beta^{\frac{\log j}{\log \beta_c}+1} + \beta -2}
\right) \\
& = &
\left(\frac{\beta-1}{2(3+\beta)}\right)^{-1}
\left(\frac{\beta-1}{j^{\frac{\log \beta}{\log \beta_c}+1} + \beta -2}
\right)
\, .
\end{eqnarray*}
The assumption $\beta>\beta_c$ makes the last expression summable
over $j$. Hence
\[
\sum_{j=2}^\infty \Pf[T_{j, tot} > T^*_{j, tot}] < \infty
\]
so that by using Borel--Cantelli yet again we get a.s.\ that
$T_{j, tot} > T^*_{j, tot}$ for at most
finitely many $j$. This takes us from the already-established
(\ref{eq:need_to_show^*}) to the desired
(\ref{eq:need_to_show_for_large_bias}), and we are done.
$\Cox$

\section{Main construction: zero speed regime}  \label{sect:zero_speed}

Having established, in the previous section, the $\beta>\beta_c$ part of
Theorem \ref{th:main}, it only remains to prove the $\beta<\beta_c$ part.

\medskip\noindent
{\bf Proof of Theorem \ref{th:main}, case $\beta<\beta_c$.}
As usual, we assume (without loss of generality due to
Lemma \ref{lem:reduction}) that the origin sits
on a branch of order $1$.
We need to show that for $\beta<\beta_c$ we have a.s.\
$\lim_{t \rightarrow \infty}\frac{X_t}{t} =0$. For this it is
enough to show that a.s.
\begin{equation}  \label{eq:takes_a_long_time}
\lim_{x \rightarrow \infty} \frac{U(x)}{x} \, = \, \infty
\end{equation}
where, as before, $U(x)$ is the time of
first arrival to $x$-coordinate $x$. To this end, we proceed as in
the last part of the proof of Proposition \ref{prop:warm-up_construction},
writing $A_k$ for the event that the first time the random walk reaches
the anchor of trap number $k$, it enters the trap and spends at least time
$\beta^{c_k}= \beta^{3^{(k-1)(k-2)/2}}$ there. Furthermore, defining
$W_i$ as the number of events happening amongst
$A_{2^{i-1}+1}, A_{2^{i-1}+2}, \ldots, A_{2^i}$
(recall (\ref{eq:def_Wi})), we get using the same estimates as those
leading up to (\ref{eq:Poisson_Wi}) that $\Pf(W_i=0)$ decays exponentially
in $i$. Hence, Borel--Cantelli tells us that a.s.
\begin{equation}  \label{eq:how_things_should_be_for_large_i}
W_i> 0 \, \mbox{ for all but finitely many } i \, .
\end{equation}
Write $x_i$ for the (random)
$x$-coordinate at which trap number $2^i$ attaches
to the path $P$, and note that $x_i$ does not exceed
$b_{2^i}= 4 \cdot 3^{2^i(2^i-1)/2}$. We have on the event
$\{W_i>0\}$ that
\[
U(x) > \beta^{3^{2^{i-1}(2^{i-1}-1)/2}} \, \mbox{ for all } \, x>x_i
\]
and consequently that
\begin{eqnarray} \nonumber
\frac{U(x)}{x} & \geq & \frac{1}{x_{i+1}} \beta^{3^{2^{i-1}(2^{i-1}-1)/2}} \\
& \geq &
\frac{\beta^{3^{2^{i-1}(2^{i-1}-1)/2}}}{4 \cdot 3^{2^{i+1}(2^{i+1}-1)/2}}
\label{eq:triply_vs_doubly_exponential}
\end{eqnarray}
for all $x \in [x_i, x_{i+1}]$. This bound tends to $\infty$ as
$i \rightarrow \infty$. Using (\ref{eq:how_things_should_be_for_large_i}),
we thus get (\ref{eq:takes_a_long_time}), so the proof is complete.
$\Cox$

\bibliographystyle{alea2}
\bibliography{07-02}

\end{document}